\documentclass[11pt]{article}

\usepackage{epsfig,amsmath,latexsym}
\usepackage{amsfonts}

\newtheorem{theorem}{Theorem}[section]

\newtheorem{prop}[theorem]{Proposition}

\def\slfrac#1#2{\hbox{\kern.1em %
 \raise.5ex\hbox{\the\scriptfont0 #1}\kern-.11em %
 /\kern-.15em\lower.25ex\hbox{\the\scriptfont0 #2}}}

\newcommand{\beq}{\begin{eqnarray}}
\newcommand{\eeq}{\end{eqnarray}}
\newcommand{\beql}[1]{\begin{eqnarray}\label{#1}}
\newcommand{\beqs}{\begin{eqnarray*}}
\newcommand{\eeqs}{\end{eqnarray*}}
\newcommand{\eqn}[1]{(\ref{#1})}

\newcommand{\la}{\lambda}
\newcommand{\af}{\alpha}

\newcommand{\rr}{{\mathbb R}}

\newcommand{\bv}{{\mathbf v}}
\newcommand{\bw}{{\mathbf w}}

\newcommand{\sA}{{\mathcal A}}

\newcommand{\sD}{{\mathcal D}}

\newcommand{\sI}{{\mathcal I}}
\newcommand{\sK}{{\mathcal K}}
\newcommand{\sT}{{\mathcal T}}

\newcommand{\sP}{{\mathcal P}}

\newcommand{\bsq}{\vrule height .9ex width .8ex depth -.1ex}

\makeatletter
\def\@sect#1#2#3#4#5#6[#7]#8{\ifnum #2>\c@secnumdepth
     \def\@svsec{}\else
     \refstepcounter{#1}\edef\@svsec{\csname the#1\endcsname.\hskip .75em }\fi
     \@tempskipa #5\relax
      \ifdim \@tempskipa>\z@
        \begingroup #6\relax
          \@hangfrom{\hskip #3\relax\@svsec}{\interlinepenalty \@M #8\par}%
        \endgroup
       \csname #1mark\endcsname{#7}\addcontentsline
         {toc}{#1}{\ifnum #2>\c@secnumdepth \else
                      \protect\numberline{\csname the#1\endcsname}\fi
                    #7}\else
        \def\@svsechd{#6\hskip #3\@svsec #8\csname #1mark\endcsname
                      {#7}\addcontentsline
                           {toc}{#1}{\ifnum #2>\c@secnumdepth \else
                             \protect\numberline{\csname the#1\endcsname}\fi
                       #7}}\fi
     \@xsect{#5}}
\def\@begintheorem#1#2{\it \trivlist \item[\hskip \labelsep{\bf #1\ #2.}]}

\def\plain{plain}\ifx\fmtname\plain\csname fi\endcsname
     
     \input docstrip
     \preamble

     Do not distribute the stripped version of this file.
     The checksum in the header refers to the documented version.

     \endpreamble
     \generateFile{here.sty}{t}{\from{here.doc}{}}
     \endinput
\fi
\ifcat a\noexpand @\let\next\relax\else\def\next{%
    \documentstyle[here,doc]{article}\MakePercentIgnore}\fi\next
\ifx\@Hxfloat\@Hundef\else\expandafter\endinput\fi
\let\@Hxfloat\@xfloat
\def\@xfloat#1[{\@ifnextchar{H}{\@HHfloat{#1}[}{\@Hxfloat{#1}[}}
\def\@HHfloat#1[H]{%
\expandafter\let\csname end#1\endcsname\end@Hfloat
\vskip\intextsep\vbox\bgroup\def\@captype{#1}\parindent\z@
\ignorespaces}
\def\end@Hfloat{\egroup\vskip \intextsep}

\makeatother

\catcode`\@=11

\catcode`@=12

\begin{document}


\begin{center}
{\Large {\bf Bounds for the $3x+1$ Problem using Difference Inequalities }}\\
\vspace{1.5\baselineskip}

{\em Ilia Krasikov}\\
\vspace*{.2\baselineskip}
Brunel University \\
Department of Mathematical Sciences \\
Uxbridge UB8 3PH, United Kingdom \\
\vspace{1.5\baselineskip}

{\em Jeffrey C. Lagarias } \\
\vspace*{.2\baselineskip}
AT\&T Labs  - Research ,  \\
Florham Park, NJ 07932-0971 \\
\vspace*{1.5\baselineskip}

(April 30, 2002) \\
\vspace*{1.5\baselineskip}
{\bf ABSTRACT}
\end{center}
We study difference inequality systems for the $3x+1$ problem
introduced by the first author in 1989. These systems can be used to
derive lower bounds for the number of integers below $x$ for
which the $3x+1$ conjecture is true. Previous methods 
obtaining such lower bounds gave away some information 
in these inequalities; we give an improvement which 
(apparently) extracts full information from the inequalities 
in obtaining a lower bound.

We deduce, by computer-aided proof,  
that for any fixed positive integer $a$ not divisible by
$3$, and for large enough $x$ (depending on $a$), 
at least $x^{0.84}$ of integers below $x$
have $a$ in their forward orbit under the $3x+1$ function. 

\vspace{1.5\baselineskip}


\setlength{\baselineskip}{1.0\baselineskip}

%
%
%
%
%
%
%

\section{Introduction}
\setcounter{equation}{0}
The $3x+1$ problem concerns the iteration of the $3x+1$ function
$T(n) = n/2$ if $n$ is an even integer, $(3n+1)/2$ if $n$ is
an odd integer. The well known $3x+1$ Conjecture asserts that all integers
$n \geq 1$ eventually reach $1$ under iteration of 
the $3x+1$ function. Known results on this problem are surveyed in
Lagarias \cite{La85} and Wirsching \cite{Wi98}.

Let $\pi_1(x)$ count the number of integers below $x$ that eventually
reach $1$ under this iteration. 
There are several methods known 
for establishing  lower bounds 
of the form $\pi_1(x) > x^\beta$ for
a positive constant $\beta$. see \cite{AL95a}, \cite{AL95b}.
The first such bound was obtained
in 1978 by Crandall ~\cite{Cr78}.
The strongest of these methods at present is
one introduced by the first author in
1989 (\cite{Kr89}), which uses systems of  difference
inequalities, and here we consider it further.  

This method formulates, for each $k \ge 2$, 
a system $\sI_k$ of functional 
difference inequalities
$(\bmod~  3^k)$, containing about $3^k$ variables, 
which certain functions, computed from $3x+1$ iterates,
satisfy; they are specified in \S2. One can  
establish  an exponential lower bound for the growth rate of  positive
monotone solutions to these inequalities, and this translates
into lower bounds for $\pi_1(x)$ of the form $x^\beta$ for some positive
$\beta.$  The paper ~\cite{Kr89} used the system
$k=2$ to obtain a lower bound $x^{0.43}$ for the number of such
integers. Later Wirsching~\cite{Wi93} used the
system $k=3$ to obtain the lower bound $x^{0.48}$, for all sufficiently
large $x$.

In 1995 Applegate and 
Lagarias \cite{AL95b} gave an approach using  nonlinear programming
to systematically extract lower bounds from the  difference inequalities 
$\sI_k$. The first step was to iterate the inequalities to obtain a
derived system of difference inequalities $\sD$ such that any 
positive, monotone solution to
the original inequalities would remain a solution of the derived
inequalities. This step 
 can be done in many ways, in an ad hoc fashion.
Given such a system of difference inequalities $\sD$ they associated 
a family of auxiliary linear
programs $L_k^{\sD}(\lambda)$ depending on a parameter $\lambda$.
The parameter  $\lambda$ lies in the interval $1 \le \lambda \le 2,$ 
and the  coefficients of
the linear program depend nonlinearly
on $\lambda$. 
If the system $\sD$ contained only ``retarded'' variables (as defined
below) then any positive
feasible solution to the linear program for a fixed $\lambda$
yields  rigorous exponential lower bounds for the
growth of any positive  montone solution of the system $\sD$, 
with exponential
growth constant $\lambda$; the associated exponent in lower bounding
$\pi_1(x)$ is then $\gamma = \log_2 \lambda.$
One searches for the largest 
value of $\lambda$ for which a positive feasible solution exists,
which is a nonlinear programming problem.
To obtain an inequality system with retarded variables only,
Applegate and Lagarias \cite{AL95b} found it necessary to apply a
 ``truncation'' operation
which weakens the inequalities and presumably weakens the 
exponential lower bounds attained. 
 Using the system $k =9$, and a particular sequence of
reductions to derive  a suitable system $\sD$, a large computation  
yielded a lower bound $\pi_1(x) \ge x^{0.81}$ for all
sufficiently large $x$. Up to now this is the best asymptotic lower
bound for $\pi_1(x)$.

The nonlinear programming approach in \cite{AL95b} does not apply
directly to the original difference inequalities because they
contain terms with  ``advanced'' variables (as defined below). 
The purpose of this paper is to establish that the lower bounds
derived from the auxiliary linear program family
associated to  the original inequality
system $\sI_k$, denoted
 $L_k^{NT}(\lambda)$, do give legitimate lower bounds for
the $3x+1$ function, even though this system contains
advanced variables. The main theorem is stated in \S2.

This result yields an immediate improvement of the exponent for
lower bounds for the $3x+1$ problem for the system $k=9$, relying on 
computations reported in   \cite{AL95b}.
The computations given there for the linear program denoted 
$L_{\lambda}^{NT}$ for $k=9$ (which has  equivalent growth expoonent to 
the linear program 
$L_9^{NT}(\lambda)$ studied here) 
yield a better exponent than any of lower bounds rigorously
established  in \cite{AL95b}.
Using a further computation 
for $k=11$, we obtain  the improved lower bound
$$\pi_1(x) > x^{0.84},$$ 
valid for all sufficiently large $x$, given in \S6.

The main interest in the improved result, however, is 
that the linear program families $L_k^{NT}(\lambda)$, 
although of exponential size in $k$,  have
a relatively simple  structure.
One hopes 
that a bound of the form $\pi_1(x)> x^{1 - \epsilon}$ 
for any $\epsilon > 0$ can eventually 
be proved  by considering $L_k^{NT}(\lambda)$ for arbitrarily
large $\lambda$, and understanding better the structure of the
feasible solution sets to these linear program systems.

%
%
%
%
%
%
%
%

\section{Main Result}
\setcounter{equation}{0}

We first recall the difference inequalities $\sI_k$ of
Krasikov~\cite{Kr89}.
We consider the $3x+1$ function $T(n)$, and for $a \not\equiv 0~(\bmod~3)$
and $x \ge 1$ we define the function
$$
\pi_a (x) :=  \#\{n: 1 \le n \le x, \mbox{some}~T^{(j)}(n) = a.\}
$$
and the related function
$$
\pi_a^*(x) := \# \{n: n \le x, ~\mbox{some}~T^{(j)}(n) = a, ~\mbox{all}~
T^{(i)}(n) \le x~ \mbox{for}~ 0 \le i \le j \}.
$$
Note that $\pi_a^*(x) \le \pi_a(x).$
For each residue class $m ~(\bmod~3^k)$ with $m \ne 0~(\bmod~3)$, 
we define for $y \ge 0$ the function
$$
\phi_k^m(y) := \inf \{ \pi_a^*(2^y a):~ a \equiv m~(\bmod~3^j) ~\mbox{and}
~a~\mbox{not in a cycle} \}.
$$
This function is well defined because there always exists
some $a \equiv m~(\bmod~3^k)$ not in
a cycle. 

This definition immediately implies that for $k \ge 2$
and all $m ~(\bmod~3^k)$, $m \not\equiv 0~(\bmod~3)$, 
these functions satisfy the three properties:

(P1) ({\em Positivity}) {\em For all} $y \ge 0$,
$$ \phi_k^m(y) \ge 1.$$

(P2) ({\em Monotonicity}) {\em For} $y \ge 0$,
$$ 
\phi_k^m(y) ~\mbox{is a nondecreasing function of}~ y.
$$

(P3) ({\em Minimization}) {\em For} $m \in [3^{k-1}]$
{\em and all} $y \ge 0$, 
$$
\phi_{k-1}^m(y) = \min[ \phi_k^m(y), \phi_k^{m + 3^{k - 1}}(y),
\phi_k^{m + 2\cdot 3^{k - 1}}(y)].
$$

It is easy to see that
\beql{200}
\phi_k^m(y) = \phi_k^{2m}(y - 1)~~\mbox{if}~ m \equiv 1~ (\bmod~ 3),
\eeq
hence it suffices to study $\phi_k^m(y)$ for $y \equiv ~2(\bmod~3).$
For convenience in what follows we  let $[3^k]$ denote the set
of congruence classes
\beql{200b}
[3^k] := \{ m~(\bmod~3^k):~ m \equiv 2 ~(\bmod~3) \}.
\eeq

The difference inequality system of 
Krasikov~\cite{Kr89} can be put in
the following form.

\begin{prop}~\label{pr21}
Let $\alpha := \log_2 3 \simeq 1.585$. For each $k \ge 2$,
the set of functions 
$\{\phi_k^m(y): m \in [3^k]\}$ 
satisfy the following system $\sI_k$ of difference inequalities,
valid for all $y \ge 2$.

 (D1) If $m \equiv 2~(\bmod~9)$ then
\beql{eq201}
\phi_k^m(y) \ge \phi_k^{4m}(y-2) + \phi_{k-1}^{(4m - 2)/3}(y+\alpha -2).
\eeq

(D2) If $m \equiv 5~(\bmod~9)$ then
\beql{eq202}
\phi_k^m(y) \ge \phi_k^{4m}(y-2).
\eeq

(D3) If $m \equiv 8~(\bmod~9)$ then
\beql{eq203}
\phi_k^m(y) \ge \phi_k^{4m}(y-2) + \phi_{k-1}^{(2m - 1)/3}(y+\alpha -1)
\eeq

In these inequalities the functions $\phi_{k-1}^m (y)$ are defined by
\begin{equation}\label{eq204}
\phi_{k-1}^m(y) := \min[ \phi_k^m(y), \phi_k^{m + 3^{k - 1}}(y),
\phi_k^{m + 2\cdot 3^{k - 1}}(y)].
\end{equation}
\end{prop}

\noindent\paragraph{Proof.}
This follows from \cite[Lemma 4]{Kr89}, 
and appears in \cite[Prop. 2.1]{AL95b}.~~~$\bsq$

We regard the system $\sI_k$ of inequalities as expressed
entirely in terms of the functions
$\{\phi_k^m (y): m \in [3^k ]\}$, 
by using the minimum formulas \eqn{eq204}.
In that case all functions appearing are of the form 
$\phi_k^m (y+ \beta_j )$ for various real numbers $\beta_j$.
If $\beta_j \ge 0$ we call such a term {\em advanced}, 
while if $\beta < 0$ we call such a term {\em retarded},
since the terms have advanced arguments and retarded arguments
respectively, in terms of the ``time'' variable  $y$.

Applegate and Lagarias~\cite{AL95b} associated to $\sI_k$ various
auxiliary linear programs $L_k^\sD(\lambda)$ depending
on a parameter $\lambda > 1$; strictly positive feasible solutions for
admissible  linear programs for a given $\lambda$ lead to exponential
lower bounds for the functions $\phi_k^m(y) \ge c_0 \lambda^y$.
Here we study  a particular linear program family,
denoted $L_\lambda^{NT}$ in \cite{AL95b}, associated
directly to the inequalities $\sI_k$, to which their
methods did not apply to get any lower bound.
Actually we consider a 
modified linear program family $L_k^{NT}(\lambda)$ given 
below; this differs from $L_\lambda^{NT}$ 
in having a different objective function variable, being a minimization
rather than a maximization, and having the certain 
nonnegativity constraints modified to make them strictly positive.
However, as shown in \S6,
this modified linear program is  equivalent to 
the one in \cite{AL95b} in the sense that matters here: it 
has a feasible solution for $\lambda$
if and only if the
corresponding linear program in \cite{AL95b} has a  strictly positive
feasible solution for the same $\lambda$. 

The linear program family  $L_k^{NT}(\lambda)$
is as follows.

\begin{equation}\label{eq205}
L_k^{NT} (\lambda): ~\mbox{Minimize }~ C_k^{max}
\end{equation}
subject to:

\noindent
(L0) {\em For all} $m \in [3^k]$,
\beql{eq206}
1  \le  c_k^m \le C_k^{max}
\eeq

\noindent
(L1) {\em For all $m \in [3^k]$ with } $m \equiv 2$ $(\bmod~ 9)$,
\begin{equation}\label{eq208}
c_k^m \le c_k^{4m} \la^{-2} + c_{k-1}^{\frac{4m-2}{3}} \lambda^{\alpha -2} ~.
\end{equation}

\noindent
(L2) {\em For all $m \in [3^k]$ with}  $m \equiv 5$ $(\bmod~9)$,
\begin{equation}\label{eq209}
c_k^m \le c_k^{4m} \la^{-2}~.
\end{equation}

\noindent
(L3)
{\em For all $m \in [3^k]$ with} $m \equiv 8$ $(\bmod~9)$,
\begin{equation}\label{eq210}
c_k^m \le c_k^{4m} \lambda^{-2} + c_{k-1}^{\frac{2m-1}{3}} \la^{\alpha -1} ~.
\end{equation}

\noindent
(L4) {\em For all} $m \in [3^k]$,
\begin{eqnarray}\label{eq211}
 c_{k-1}^m & \le & c_k^m ~, \\
\label{eq212}
 c_{k-1}^m & \le & c_k^{m+3^k} ~, \\
\label{eq213}
 c_{k-1}^m & \le & c_k^{m+2 \cdot 3^k} ~.
\end{eqnarray}

Note that the inequality signs in (L1)--(L3) go in the 
opposite direction from that in
the difference inequalities $\sI_k$,
while (L4) goes in the same direction.

We call the  variables $\{c_k^m: m \in [3^k ]\}$ 
{\em principal variables} in $L_k^{NT} (\la )$,
and  the variables $\{c_{k-1}^m : m \in [3^{k-1}] \}$ 
{\em auxiliary variables}; 
the remaining variable  $C_k^{max}$ is the objective
function variable. The objective function
variable itself plays no role in determining feasibility of
the linear program; the inequalities it appears in
can always be satisfied by
setting it equal to the maximum of the principal variables.
If this linear program has any feasible solution,
then this solution may be rescaled by a multiplicative constant 
so that $\min{ \{c_k^m \} } = 1$, while decreasing $C_{max}$,
hence any  optimum value of this linear program will have  
$\min{ \{c_k^m \} } = 1$. Given a feasible solution, set
\beql{eq214}
\bar{c}_{k-1}^m := 
\min \{ c_k^m,~ c_k^{m + 3^{k-1}},~ c_k^{m + 2\cdot 3^{k-1}} \}.
\eeq
The inequalities (D4) say that $c_{k-1}^m \le \bar{c}_{k-1}^m.$
There are no lower bounds imposed on the auxiliary variables
$c_{k-1}^m$, but given any feasible solution, there exists a 
positive
feasible solution with the same principal variables and with
auxiliary variables
$$ c_{k-1}^m = \bar{c}_{k-1}^m \ge 1.$$
Indeed (D4) still holds for this choice of auxiliary variables
and the remaining inequalities (D1)-(D3) stay the same or weaken.

The linear program $L_{k}^{NT}(\lambda)$
encodes advanced variables, and the theorems in
\cite{AL95b} do not apply to it.
Conjecture 4.1 of \cite{AL95b} asserts that the
largest value of $\lambda$ for which
$L_{k}^{NT}(\lambda)$ has a positive feasible solution should
give the largest possible exponential lower bound
for positive, monotone functions $\Phi_k$ satisfying $\sI_k$.
Our main result is that $L_{k}^{NT}(\lambda)$     gives legitimate
lower bounds for positive solutions 
for such functions  $\phi_k^m(y)$. 

\begin{theorem}~\label{th21}
Let $1 \le \lambda \le 2$ be such that the linear program
$L_k^{NT}(\lambda)$ has a feasible solution
with principal variables $\{c_k^m : m \in [3^k]\}$. 
Then for and all $m \in [3^k]$ and all $y \ge 0,$ 
\beql{eq217}
\phi_k^m(y) \ge \Delta_1 \cdot c_k^m \lambda^y,
\eeq
in which
\beql{eq218}
\Delta_1 := \frac{1} 
{4 \max{ \{c_k^m: m \in [3^k]\}} }.
\eeq
\end{theorem}

We believe, although we have no proof, that this result
gives the largest exponential-type lower bound that
can be extracted from the difference inequalities $\sI_k$.
This is discussed at the end of \S6.

Theorem~\ref{th21} is established as follows. In \S3 we show that there
exists a sequence of back substitutions of the difference inequalities
into themselves that results in a difference inequality system
from which all advanced variables have been eliminated. This
results in a new system of difference inequalities $\sI_k(EL)$.
We show that all solutions $\phi_k^m$ of $\sI_k$
which  possess the positivity and monotonicity properties (P1) and (P2)
will also be solutions of  $\sI_k(EL)$.

In \S4 we consider linear programs.
To each difference inequality system $\sD$ (of a specified kind)
we associate in
a strictly deterministic way an 
auxiliary linear program family $L^{\sD}(\lambda)$. Let
$L_k^{EL}(\lambda)$ denote  the linear program family attached 
to $\sI_k(EL)$. 
The main result of \S4 is the deduction that
if the linear program $L_k^{NT}(\lambda)$
has a positive feasible solution with principal variables
$\{c_k^m : m \in [3^k]\}$, then 
the linear program $L_k^{EL}(\lambda)$ with the same value of
$\lambda$ also has a positive feasible solution with the
same principal variable values.

In \S5 we show that any difference inequality system $\sD$
in which only retarded variables appear has the property that
positive feasible solutions to the auxiliary linear program
$L^{\sD}(\lambda)$ for fixed $\lambda$ yields lower
bounds of the form \eqn{eq217};
the proof is similar to \cite[Theorem 2.1]{AL95b}.
 It immediately follows that we get such
lower bounds from the linear program family $L_k^{EL}(\lambda)$.
We then prove Theorem~\ref{th21}, by combining  this result
with the main result of \S4.

In \S6 we present taxonomic 
 data on the derived systems
 $L_k^{EL}(\lambda)$
for $2 \le k \le 5$ and information on positive feasible
solutions the system
 $L_k^{NT}(\lambda)$
for $2 \le k \le 11$, computed by David Applegate, which yield
the lower bound $\pi_a(x) \ge x^{0.84}$ for all sufficiently
large $x$. The results of \S4 imply that the linear program
family  $L_k^{EL}(\lambda)$  might  conceivably give better
exponential lower bounds than are obtainable from the 
linear program family  $L_k^{NT}(\lambda)$. Numerical
experiments show this is not the case for $2 \le k \le 5$; 
here $k=5$ was the limit of computability for the
system $L_k^{EL}(\lambda)$.

%
%
%
%
%
%
%
%

\section{Eliminating Advanced Variables}
\setcounter{equation}{0}
We describe a recursive back-substitution procedure to eliminate
``advanced'' terms of the inequality system $\sI_k$.
We view the inequality system $\sI_k$ as expressed 
entirely in terms of functions $\phi_k^m (y+\beta )$ 
 by replacing each term 
involving any variable $\phi_{k-1}^{m'}(y + \beta')$
by the  minimization expression 
on the right side of \eqn{eq204}
in terms of $\phi_k^m$ functions.

We start with a single inequality (D3) of the system $\sI_k$ 
associated to a fixed $m \in [3^k]$, $m \equiv 8$ $(\bmod~9)$, 
and perform a recursive back-substitution process 
of the inequalities $\sI_k$ into its right-hand side.
At the $l$th-stage of this process we will have 
an inequality $I_k^m (l)$ whose left side is $\phi_k^m (y)$ 
and whose right side is a nested series
of minimizations of various functions $\phi_k^{m'} (y+ \beta' )$.
The step from $I_k^m (l)$ to $I_k^m (l+1)$ has two substeps.
First, one picks
 an advanced term $\phi_k^{m'} (y+ \beta' )$,
$\beta' \ge 0$ appearing in the right side of $I_k^m (l)$ 
and replaces it with the right side of the inequality 
$\sK_k$ for $\phi_k^{m'} (y')$ with $y' = y+ \beta'$.
(This is called ``splitting'' a term in \cite{AL95b}.)
A new minimization term may appear in this process, 
which contains three terms
\beql{eq31}
\phi_k^{m'} (y+ \beta'') , \phi_k^{m'+3^{k-1}} (y+\beta''),
\phi_k^{m' + 3^{k-1}} (y+ \beta''' ) \,.
\eeq
The second substep in obtaining $I_k^m (l+1)$ is to apply a 
deletion rule described below, which, if $\beta'' \ge 0$, 
may remove up to two of these terms. The resulting inequality
after the deletion substep is $I_k^m (l+1)$.

At each stage in this process the inequality $I_k^m (l)$ has 
$\phi_k^m (y)$ on its left side and a sum of nested minimization 
terms on its right side, involving various functions $\phi_k^m (y+ \beta_j)$;
it will have each  $\beta_j \ge -2$, because we will only 
substitute for terms $\phi_k^m (y+ \beta_j)$ with 
$\beta_j \ge 0$, and the formulas (D1)--(D3) produce new terms
$\phi_k^{m'} (y+ \beta'_j)$ which have $\beta'_j \ge \beta_j -2$.
The structure of the right side of an inequality
$I_k^m (l)$ is described by a directed rooted labelled 
tree $\sT_k^m (l)$, in which the root mode is labelled 
with the left side $\phi_k^m (y)$ of the original inequality, 
each node is either a $p$-node (for ``principal'') 
or an $m$-node (for ``minimization'').
The initial tree for the  inequality $I_k^m$ for
an $m \in [3^k]$ with $m \equiv~8~(\bmod~9)$ is pictured in 
Figure 1.

%
%

\begin{figure}[htb]
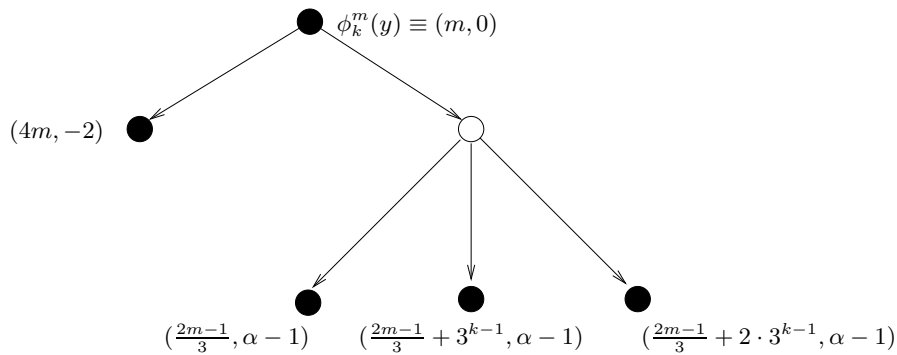

\begin{center}
\input kafg3.pstex_t
\end{center}
\caption{Rooted tree for inequality (D3).}
\end{figure}~\label{figure1}

Here $p$-nodes are indicated by solid points and $m$-nodes
by circled points.
Each $p$-node is labelled by data $(m', \beta')$ specifying 
the function $\phi_k^{m'} (y+ \beta' )$, while each $m$-node 
is labelled by data $(m',\beta')$ of the
$p$-node of which it is a child.
The root node is a $p$-node and has  label $(m,0)$.
The inequality $I_k^m (l)$ is uniquely specified by the tree
$\sT_k^m (l)$ and vice-versa;
the root node specifies the left side $\phi_k^m (y)$ of 
the inequality $I_k^m (l)$, leaf nodes specify functions 
appearing in the right side, and the internal tree structure 
specifies the nested sequence of additions and minimizations 
comprising this right side of the inequality.

A step from $\sT_k^m (l)$ to $\sT_k^m (l+1)$ consists of 
picking a leaf node with label $\phi_k^{m'} (y+ \beta' )$ 
which has $\beta ' \ge 0$ and changing the
tree in the following two substeps. First we  attach to thie leaf node 
(as root node) the directed tree associated to the formula
(D1)--(D3) of $\phi_k^{m'} (y')$ with $y' = y+ \beta'$.
We term this ``splitting'' the leaf node, following \cite{AL95b}.
The tree $\tilde{\sT}_k^m (l+1)$ that results has a new $p$-node labelled
$\phi_k^{4m'} (y+ \beta' -2)$, and may or may not have 
a new $m$-node with three new leaves \eqn{eq31} depending from it.
If there is no $m$-node this tree will be  $\sT_k^m (l+1)$.
Second, if there is a new $m$ node, we apply the {\em deletion rule} 
given below to $\tilde{\sT}_k^m (l+1)$ 
to remove  some
(possibly empty) subset of the 
three leaves in a $m$-term. 

\noindent{\bf Deletion Rule.} {\em 
For each such leaf $\phi_k^{m''} (y+ \beta'' )$, 
if $\beta'' \ge 0$, consider the directed path 
from the root node $\phi_k^m (y)$ to this leaf.
At each internal $p$-vertex on this path, one has an 
associated value $\phi_k^{m'} (y+ \beta' )$.
The leaf node $\phi_k^{m''} (y+ \beta'' )$ is deleted 
if there is an internal node with $m' \equiv m''$ $(\bmod~3^k )$ 
and with $\beta < \beta''$.}

After the deletion rule is applied to $\tilde{\sT}_k^m (l+1)$, 
the tree that results is $\sT_k^m (l+1)$, and
the inequality correspoding to it is $I_k^m (l+1)$.
We will show that the deletion rule cannot remove 
all three leaves, hence all leaf nodes on the new tree 
$\sT_k^m (l+1)$ are $p$-nodes, so the process can continue.

The back-substitution process is not completely specified, 
in that one has the
freedom to choose to split any leaf node carrying an advanced term.
In practice it is convenient to require that all nodes at a given depth 
that have advanced terms be split before proceeding 
to split nodes at greater depth, but the order of splitting
does not matter as the following result asserts.

\begin{theorem}\label{th31}
Let $k \ge 2$, and take $m \in [3^k]$ with $m \equiv 8$ $(\bmod~9)$.
The back-substitution process applied to $\phi_k^m (y)$ 
halts after a finite number of steps at an inequality 
$I_k^m (l)$ having no advanced terms on its right side.
The number of steps $l$ and the final inequality $I_k^m (l)$ 
are independent of the order in which advanced terms are split;
let $I_k^m (EL)$ denote this final inequality.
\end{theorem}

\noindent{\bf Proof.}
We first show that 
 the back-substitution procedure always halts.
We suppose not, and obtain a contradiction.
Let $\sT_l \equiv \sT_k^m(l)$ denote the rooted labelled tree associated to 
the inequality $I_k^m (l)$ for $l=1,2,\ldots$.
Then we have an infinite sequence of trees, 
each containing the last as a subtree having the same root, 
and the process defines an infinite 
limiting tree $\sT_\infty$. 
Without loss of generality we  can suppose that $\sT_\infty$ 
 has the property that in it all nodes that 
can be split are split, if necessary by doing additional splittings
of any advanced nodes that were missed, using transfinite induction.
By Konig's infinity lemma there is an infinite directed path 
in  $\sT_\infty$ starting from the root.
Along that path there is some residue class $m' \in [3^k]$ 
that occurs as a label infinitely often.
Let $\phi_k^{m'} (x+ \beta_j )$ be the successive labels 
of the $p$-nodes on this path having residue class 
$m'$ $(\bmod~3^k )$, starting from the root.
We must have each $\beta_j \ge 0$ (or the process halts) and also
\beql{eq32}
\beta_1 > \beta_2 > \beta_3 > \cdots ~,
\eeq
because the deletion rule would have removed the $p$-node 
labelled $\phi_k^{m'} (x+ \beta_j )$ if $\beta_j \ge \beta_i$ 
for some $j > i$.

The tree $\sT_\infty$
has a recursive self-similar structure, using the fact that
all nodes that could be split were split.
Consider the subtree $\sT_\infty [j]$ grown starting from 
the root node $\phi_k^{m'} (y+ \beta_j )$ along this chain, 
using the variable
$y_j = y+ \beta_j$.
These subtrees are all identical, and $\sT_\infty [2]$ is 
obtained from $\sT_\infty [1]$ by shifting the argument of $y$
by $\delta = \beta_2 - \beta_1 > 0$.
The isomorphism of $\sT_\infty [2]$ and $\sT_\infty [1]$ 
identifies $\sT_\infty [j]$ with $\sT_\infty [j-1]$, and therefore,
by induction on $j \ge 2$, we obtain $\beta_j - \beta_{j-1} = \delta$.
Thus $\beta_j = \beta_1 + (j-1) \delta$ for all $j \ge 2$, hence
$\beta_j < 0$ for sufficiently large $j$,
which contradicts all $\beta_j \ge 0$.

The back-substitution process halts at a unique tree, 
regardless of the order leaf nodes are split, 
because the back-substitution process on a given leaf
node $\bv$ does not depend on any other leaf nodes, 
but only on the path from the root node to $\bv$.
One grows out all leaf nodes until they halt, 
and the total number of steps $l$ until halting is independent 
of the order of growth.~~~$\bsq$

\begin{theorem}\label{th32}
Let $\sI_k(EL)$ denote the difference inequality system 
consisting of the inequalities {\rm (D1), (D2)} of $\sI_k$ 
plus the complete set of inequalities
$\{I_k^m (EL): m \in [3^k]$, $m \equiv 8$ $(\bmod~9) \}$.
If $\Phi_k= \{\phi_k^m (y) : m \in [3^k]\}$
 is any set of functions in which each 
$\phi_k^m (y)$ is strictly positive and nondecreasing on 
$\rr_{\ge 0}$ and satisfies the inequality system $\sI_k$ 
for all $y \ge 2$, then $\Phi_k$ also satisfies the inequalities 
$\sI_k(EL)$ for all $y \ge 2$.
\end{theorem}

\noindent{\bf Proof.}
It suffices to show that if the set
$\Phi_k :=\{\phi_k^m (y): m \in [3^k]\}$
of positive nondecreasing functions on $\rr^+ = \{ y \ge 0 \}$ 
satisfies $\sI_k$ for all $y \ge 2$, then they satisfy each 
inequality $I_k^m (l)$ for each $l \ge 1$, for all $y \ge 2$.

We prove, by induction on $l \ge 1$, that the set $\Phi_k$ 
satisfies $\sT_k^m (l)$.
The base case $l=1$ holds because $\sT_k^m (l)$ has only one 
internal $p$-node, its root node, and the corresponding 
inequality $I_k^m (1)$ is a member of $\sI_k$.
Now suppose the induction hypothesis holds for $\sT_k^m (l)$, and consider
$\sT_k^m (l+1)$.
To obtain $\sT_k^m (l+1)$ we first split a leaf of $\sT_k^m (l)$ to obtain
a tree $\tilde{\sT}_k^m (l+1)$ and then,
if a new $m$-node was added, we apply the deletion rule to the 
three vertices of that $m$-node.
The splitting procedure yielding $\tilde{\sT}_k^m (l+1)$ 
substitutes an inequality of $\sI_k$, hence $\Phi_k$ 
automatically satisfies $\tilde{\sT}_k^m (l+1)$.

Consider the deletion step, applied to the three leaf-node 
labels of $\tilde{\sT}_k^m (l+1)$ inside the min-term
\beql{eq304}
f(y) := \min [ \phi_k^{m'} (y+ \beta') , \phi_k^{m''} (y+ \beta' ) , 
\phi^{m'''} (y+ \beta' ) ] ~,
\eeq
in which $m' \equiv m'' \equiv m'''$ $(\bmod~3^{k-1} )$.
The leaf node $\phi_k^{m''} (y+ \beta' )$ is to be deleted 
if earlier in its directed path from the root appears a $p$-node 
with label $\phi_k^{m'} (y+ \beta'' )$ with $\beta < \beta'$.

To justify the deletion rule, note that the inequality 
associated to each tree $\sT_k^m (l)$ for fixed functions 
$\Phi_k$ and a fixed value $y \ge 2$, can be written as 
a sum of terms corresponding to a subset of leaves of the tree 
which are specified by choosing one of the terms 
in each min-term that attains the minimum.
(This choice is usually unique once the functions $\Phi_k$ 
and the value $y$ are specified, unless two terms 
in a min-term have equal values.)
We call this set of leaves a {\em critical assignment}, 
the leaves in it
{\em critical leaves}, and the set of paths to 
these leaves {\em critical paths}.

\noindent{\bf Claim.} {\em 
To each internal $p$-vertex $\bv$ of the tree with label 
$\phi_k^{m''} (y+ \beta'' )$, and for each fixed value of
$y \ge 2$,
one of two possibilities occurs.

(a) There are no critical assignments $\sA$ having 
a critical path passing through $\bv$.

(b) There is at least one critical assignment $\sA$ 
with a path passing through $\bv$.
For any such assignment
\beql{eq305}
\phi_k^{m''} (y+ \beta '') \ge 
\sum_{\bw \in \sA_\bv} \phi_k^{m(\bw )} (y+ \beta (\bw )) ~,
\eeq
where $\sA_{\bv}$ denotes the set of critical leaves in 
$\sA$ whose paths pass through $\bv$.
} \\

Warning: which case (a) or (b) occurs depends on the
value of $y$.
The key content of the claim is the property \eqn{eq305}
enforced in case (b). 

We will prove the claim by induction on $l$,
and justify the deletion rule at the same time.
Now \eqn{eq305} holds for the base case $l=1$ 
where the only internal $p$-node is the root node, 
and \eqn{eq305} is then an inequality in $\sI_k$.
We assume it holds for $\sT_k^m (l)$ 
and wish to prove it for $\sT_k^m (l+1)$.
First of all, the relations (a), (b) hold 
for all $\tilde{\sT}_k^n (l+1)$.
They hold for internal $p$-nodes inherited from 
$\sT_k^{(m)} (l)$, because
we have back-substituted $\sI_k$ on the right side of \eqn{eq305}.
We have added one new internal node $\bv'$, 
the one that was split, and for it 
condition \eqn{eq305} in (b) directly expresses 
the $\sI_k$ inequality substituted.

We call a vertex $\bv$ of $\tilde{\sT}_k^m (l+1)$ 
{\em totally non-critical} if no critical path passes through it,
for any critical assignment $\sA$, for any $y \ge 2$;
that is, case (a) holds for $\bv$ for all $y \ge 2$.
We can safely delete all totally non-critical vertices in 
$\tilde{\sT}_k^m (l+1)$, 
and property (b) will still hold for the resulting tree $\sT'$.
(The property that a vertex in a tree is totally non-critical is
hereditary in the sense that all vertices below a totally non-critical vertex 
are also totally non-critical.)

We now show that, for those sets of functions $\Phi_k$ 
that are positive and monotone,
all vertices removed by the deletion rule 
are totally non-critical.
Suppose the deletion rule appears to the leaf vertex $\bw$
with label $\phi_k^{m'} (x+ \beta' )$ of $\tilde{\sT}_k^m (l+1)$, 
and let $\bv$ be a vertex on its directed path that has label 
$\phi_k^{m'} (y+ \beta'' )$ with $\beta'' \le \beta'$.
If $\bw$ is not totally non-critical, there is some $y \ge 2$ 
and a critical assignment $\sA$ containing $\bw$ as a critical vertex.
Formula \eqn{eq305} of (b) applies to gives 
$\phi_k^m (y+ \beta'' ) \ge \sum_{(\tilde{m}, \tilde{\beta}) \in 
\sA_\bv} \phi_k^{\tilde{m}} (y+ \tilde{\beta} )$.
We deduce
\beql{eq306}
\phi_k^{m'} (y+ \beta'' ) \ge \phi_k^{m'} (y+ \beta' ) ~,
\eeq
because $\phi_k^{m'} (y+ \beta' )$ is the contribution of $\bw \in \sA_\bv$.
There is at least one more critical path in the sum $\sA_{\bv}$ 
which passes through the $p$-vertex $\bv'$ that was split, 
whose label is $\phi_k^{m''} (y+ \beta' )$, to its 
direct $p$-node descendant $\phi_k^{4m''} (y+ \beta' -2 )$.
Now $\beta \ge 0$ since the node $\bv'$ is split, hence $y+ \beta' -2 \ge 0$, 
so $\phi_k^{4m''} (y+ \beta' -2 ) > 0$ by positivity and 
monotonicity of $\Phi_k$.
We conclude that \eqn{eq306} can be replaced by strict inequality
\beql{eq307}
\phi_k^m (y+ \beta'' ) > \phi_k^{m'} (y+ \beta' ) ~.
\eeq
Since $\beta '' < \beta'$, this violates monotonicity of $\Phi_k$,
 the desired contradiction.

Thus, the vertices removed by the deletion rule are 
totally non-critical, hence for the resulting tree $\sT_k^m (l+1)$, 
the criteria (a), (b) and \eqn{eq305} hold for all
$p$-vertices, for the functions $\Phi_k$, for all $y \ge 2$.
This completes the claim's  induction step, and proves the claim.

Now we may apply
 \eqn{eq305} to the root vertex $\bv$ 
for all critical assignments $\sA$ for all $y \ge 2$ 
is equivalent to saying that the $\Phi_k$ satisfy 
the inequality $I_k^m (l+1)$ associated to $\sT_k^m (l+1)$ 
for all $y \ge 2$.
This completes the main induction step.~~~$\bsq$

\noindent\paragraph{Remark.}
The inequality system $\sI_k(EL)$ involves nested minimization 
to a depth $d(k)$ which grows exponentially with $k$.
The exponential growth occurs because the deletion rule requires labels 
$\phi_k^{m'} (y+ \beta_j)$,
$\phi_k^{m''} (y+ \beta_j )$ with $m' \equiv m''$ 
$(\bmod~3^k )$ and these are typically separated by distance 
comparable to $3^k$.
We present statistics in Table~\ref{table1} 
on the size of this inequality system $\sI_k(EL)$ for $2 \le k \le 5$,
computed by D. Applegate.
We measure the size in
two ways: the depth of nested minimizations, and the total of the 
number of terms that appear in such an inequality. The data is for
the term $\phi_k^m(y)$ that had the largest expansion under the
elimination procedure. 

%
%
\begin{table}[hp]\centering
\begin{tabular}{|r|r|r|}
\hline
\multicolumn{1}{|c|}{$k$} &\multicolumn{1}{c|}{\mbox{depth}}&
\multicolumn{1}{c|}{\# (\mbox{literals})} 
\\ \hline
2 & 3 & 8 \\
3 & 10 & 84 \\
4 & 41 & 12829 \\
5 & $> 226$ & $ > 10^{9}$ \\
\hline
\end{tabular}
\caption{Statistics on  $\sI_k(EL)$ Inequalities }~\label{table1}
\end{table}

%
%
%
%
%
%
%
%

\section{Linear Programs}
\setcounter{equation}{0}
We associate to a general difference inequality system $\sD_k$
(of a sort described below) 
a family of linear program $L_k^{\sD} (\lambda )$, as follows.
We suppose that $\sD_k$ consists of inequalities 
$\{D_k^m : [m] \in 3^k \}$ 
in which each inequality $D_k^m$ 
is described by a rooted labelled tree $\sT_k^m$ 
of the type considered in \S3, involving variables
$\{c_k^m: [m] \in 3^k\}$.
The linear program has the basic form:
\beql{eq401}
L_k^{\sD} (\lambda) : ~ \mbox{Minimize}~C_{max}
\eeq
subject to, for all $m \in [3^k]$,
$$1 \le c_k^m \le C_k^{max},$$
together with all inequalities associated to each tree $\sT_k^m$
as specified below.

The LP-inequality system associated to a given tree $\sT$ 
involves the  principal variables 
$\{c_k^m : m \in [3^k] \}$
and certain auxiliary variables
$\{ a_\bv : \bv ~\mbox{an $m$-vertex of $\sT$} \} ~.$
These auxiliary variables are distinct for different trees $\sT_k^m$.
We associate to each node $\bw$ the label $(m(\bw), \beta (\bw))$
which consists of a  {\em residue class} $m(\bw)$
and a {\em weight} $\beta ( \bw )$.
For a $p$-node $\bw$ these labels are 
determined by its associated function $\phi_k^{m(\bw )} (x+ \beta (\bw ))$ 
with $m ( \bw )$ determined $(\bmod~3^k )$.
For an $m$-node it is taken from the 
node function of any of its children, 
where we view $m( \bw )$ $(\bmod~3^{k-1})$ in this case, 
noting that $\beta (\bw ) (\bmod~3^{k-1})$ 
is the same for all the child nodes.
To specify the inequalities, we subdivide the tree $\sT$ into levels:
we say that a vertex $\bw$ is at {\em $m$-depth $d$} if there are exactly
$d-1$ internal $m$-nodes on the path from the root node 
to $\bw$ (not counting $\bw$ itself).
The LP inequalities associated to $\sT$ are in 
one-one correspondence with the leaf nodes of $\sT$.
To each leaf node $\bw$ we assign a rooted subtree 
$\sT_\bw$ which consists of:

(1) The terminal part of the path from the root node to the leaf node.
If an $m$-node occurs on the path, then it consists of 
that part of the path from the final $m$-node to the leaf;
if no $m$-node occurs then it is the entire path from the root.
We denote this path $\sP_\bw$ and call its top node the $\bw$-root node.
Every vertex on $\sP_\bw$ is a $p$-node except possibly the $\bw$-root node.

(2) All other children of any $p$-node on the path $\sP_\bw$.
These other children are all $m$-nodes.

\noindent
A typical subtree $\sT_\bw$ is pictured in Figure 2.

%
%

\begin{figure}[htb]
\centerline{\psfig{file=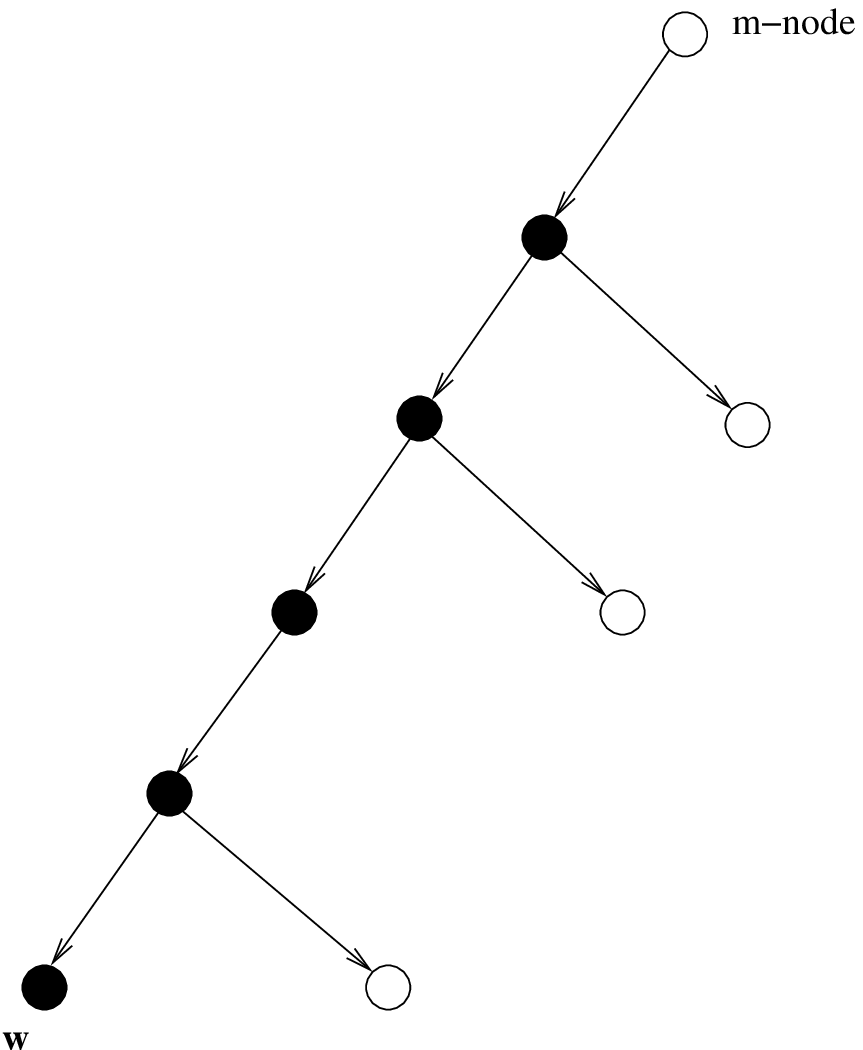,width=3in}}
\caption{Subtree $\sT_{\bw}$ of the leaf node $\bw$ 
($m$-nodes are circled)}
\end{figure}~\label{figure2}

All the edges of $\sT$ are partitioned among the $\sP_\bw$ 
and each $\sP_\bw$ contains exactly one leaf node.
The trees in this partition are also in one-one correspondence with:
either the root node $\bv$ or a pair $(\bv, \bv' )$ 
consisting of an $m$-node $\bv$ and one of its children $\bv'$.

The LP-inequality associated to the unique leaf node $\bw$ 
having no $m$-nodes on its path is of the form
\beql{eq402}
c_k^m \le \la^{\beta ( \bw )} c_k^{m(\bw )} +
\sum_{\bv \in \sT_\bw \atop \mbox{$m$-node}} \la^{\beta (\bv )}a_\bv
\eeq
where $m=m(\bv_0)$ for the root vertex $\bv_0$.
For all leaf nodes $\bw$ such that $\sT_\bw$ 
has a node $\bv_0$ as $\bw$-root node, the associated LP-inequality is
\beql{eq403}
\la^{\beta (\bv_0 )} a_{\bv_0} \le \la^{\beta (\bw )} c_k^{m(\bw )} +
\sum_{{\bv \in \sT_\bw \atop \mbox{$m$-node}} \atop \bv \neq \bv_0}
\la^{\beta (\bv )}a_\bv ~.
\eeq
Note that the direction of this LP-inequality \eqn{eq402} 
for the root node is
opposite to that of the $\phi_k^m (y)$-inequality.

For the original difference inequality system  $\sI_k$, 
the linear program $L_k^\sI (\la )$ produced in this way 
is equivalent to $L_k^{NT} (\la )$ in the following sense:
to every feasible solution of $L_k^{NT} (\la )$ with 
principal variables $\{c_k^m \}$ there corresponds 
a feasible solution to $L_k^{\sI} (\la )$ 
with the same principal variables, and vice-versa.
To see this, we note that $L_k^{NT} ( \la )$ has 
auxiliary variables $c_{k-1}^m$, while $L_k^{\sI} (\la )$ 
has auxiliary variables $a_{\bv}$ in one-one correspondence 
with $c_k^m$ for $m \equiv 2$ or 8 $(\bmod~9)$;
the LP-inequalities in $L_k^\sI (\la )$ on these variables 
are equivalent to
\beql{eq404}
a_\bv \le \bar{c}_{k-1}^{\frac{4m-2}{3}}\quad\mbox{or}\quad
a_\bv \le \bar{c}_{k-1}^{\frac{2m-1}{3}} ~,
\eeq
according as $m \equiv 2$ $(\bmod~9)$ or $m \equiv 8$ $(\bmod~9)$,
respectively, where 
$\bar{c}_{k-1}^m := \min_{0 \le j \le 2} \{ c_k^{m+j3^{k-1}} \}.$
The correspondence between feasible solutions of $L_k^{NT} (\la )$ and
$L_k^{\sI} (\la )$ is obtained by setting
\beql{eq405}
a_\bv = \bar{c}_{k-1}^{\frac{4m-2}{3}} \quad\mbox{or}\quad
a_\bv = \bar{c}_{k-1}^{\frac{2m-1}{3}} ~,
\eeq
according as $m \equiv 2$ $(\bmod~9)$ or $m \equiv 8$ $(\bmod~9)$,
respectively. 

We let $L_k^{EL} (\la )$ denote the family of linear programs 
associated to the derived
inequality system $\sI_k^m (EL )$ of Theorem \ref{th32}.

\begin{theorem}\label{th41}
Suppose for a given $\la$ with $1 \le \la \le 2$ 
that the linear program $L_k^{NT} (\la )$ has 
a feasible solution with principal variables 
$\{c_k^m: m \in [3^k]\}$.
Then the linear program $L_k^{EL} (\la )$ has 
a positive feasible solution with the same principal variables.
\end{theorem}

\noindent{\bf Proof.}
We prove this by starting with the inequality system 
$\sD_1 := \sI_k$ and then successively producing 
inequality systems
$\{\sD_j : 1 \le j \le r \}$, in which $\sD_{j+1}$ is
obtained from $\sD_j$ by a single back-substitution 
in one inequality, and ending at the
final system $\sD_r = \sI_k (EL)$.
For definiteness we choose to do 
the back-substitution procedure on each inequality $I_k^m$, 
for $m \in [3^k]\}$ in order 
until it halts, as guaranteed by Theorem \ref{th31}, 
and go to the next $m$, in the order  $m=2,5, 8, \ldots , 3^k -1$.

We prove by induction on $j \ge 1$ that if 
$\{c_k^m : m \in [3^k ]\}$ 
yields a feasible solution of $L_k^{NT} (\la )$, 
then these same principal variable values occur 
in some positive feasible solution of $L_k^{\sD_j} (\la )$.
The base case $j=1$ holds because the linear program 
$L_k^{\sD} (\la )$ agrees with $L_k^{NT} (\la )$;
when we assign the auxiliary variables $a_\bv$ 
the values \eqn{eq405} we obtain a positive feasible solution 
with the given $\{c_k^m \}$.

For the induction step, first note that in going 
from $\sD_j$ to $\sD_{j+1}$, we ``split'' one leaf vertex $\bw$ 
of a particular tree $\sT_k^m(l)$, leaving all other trees alone, 
and then perform a deletion operation.
The vertex $\bw$ being a $p$-node, 
has associated value $m( \bv )$ $(\bmod~3^k )$.
We let $\tilde{\sD}_{j+1}$ denote the inequalities 
resulting from the splitting operation before the deletion step.
It suffices to show that $L_k^{\tilde{\sD}_{j+1}}$ has a
feasible solution with the same principal variables, 
for the deletion step merely deletes 
linear programming inequalities, which preserves feasible solutions.
The splitting step changes exactly one of the inequalities in $L_k^{\sD_j}$;
if it adds a new $m$-vertex $\bv$, then it adds on up to three new
inequalities, each involving the 
new auxiliary variable $a_{\bv}$ for the added $m$-vertex.
The corresponding tree is updated to $\sT_k^m(l+1)$.

Let $m' = m ( \bw )$.
If $m' \equiv 5$ $(\bmod~9)$, the unique $LP$-inequality 
containing the term $c_k^{m'} \la^{\beta ( \bw )}$ 
corresponding to $\bw$ on its right side,
has this term replaced by that of a new leaf vertex 
$\bw'$ with $m(\bw' ) = 4m'$, $\beta(\bw') = \beta (\bw ) -2$ 
and $\bw'$ has the same depth no $\bw$;
its new term is $c_k^{4m'} \la^{\beta (\bw )-2}$.
However by hypothesis $\{c_k^m\}$ satisfies $L_k^{NT} (\la )$, 
hence it satisfies the inequality
$$c_k^{m'} \le c_k^{4m'} \la^{-2} ~.$$
Thus we obtain 
$c_k^{m'} \la^{\beta (\bw )} \le c_k^{4m'} \la^{\beta (\bw-2}$, 
so the right side of the new inequality \eqn{eq402} or \eqn{eq403} 
is less binding than before, and the solution remains feasible.
If $m' \equiv 2$ $(\bmod~9)$, the term $c_k^{m'} \la^{\beta (\bw )}$ 
is replaced with
$$c_k^{4m'} \la^{\beta (\bw )-2} + a_\bv \la^{\beta (\bw )} ~,$$
where $\beta (\bv ) = \beta (\bw ) + \alpha -2$, and 
$L_k^{\tilde{\sD}_{j+1}}$ has
three new inequalities
\beql{eq406}
a_\bv \la^{\beta (\bv )} \le c_k^{m(\bv ) + j3^{k-1}} \la^{\beta (\bv )}
\eeq
for $0 \le j \le 2$, with $m(\bv ) = \frac{4m(\bw )-2}{3}$.
We may choose
\beql{eq407}
a_\bv = \bar{c}_{k-1}^{m(\bv)} := 
\min_{0 \le j \le 2} \{c_k^{m(\bv ) + j3^{k-1}} \}
\eeq
and satisfy \eqn{eq406};
the fact that
$\{c_k^m \}$ satisfies $L_k^{NT} (\la )$ gives
$$c_k^{m'} \la^{\beta (\bw )} \le
c_k^{4m'} \la^{\beta (\bw )-2} + c_{k-1}^{\frac{4m'-1}{3}}
\la^{\beta (\bw ) + \af _2} =
c_k^{4m'} \la^{\beta (\bw )-2} - a_\bv \la^{\beta (\bv )} ~.
$$
Thus the right side of the equation is less binding than before,
so remains feasible.
The case $m' \equiv 8$ $(\bmod~9)$ is handled by similar reasoning 
to the case $m' \equiv 2$ $(\bmod~9)$, 
so feasibility is maintained in this case.
The induction step follows.

The final case of the induction step gives the 
inequality system $\sI_k (EL)$, and the theorem follows.~~~$\bsq$

%
%
%
%
%
%
%
%

\section{Lower Bounds For Difference Inequalities}
\setcounter{equation}{0}

We obtain exponential lower bounds for systems of 
positive nondecreasing functions $\Phi_k$ satisfying 
difference inequalities $\sD$ without advanced variables,
using an associated linear program $L_k^{\sD}$.
The following result is similar in spirit 
to \cite[Theorem 2.1]{AL95b}.

\begin{theorem}\label{th51}
Let $\Phi_k := \{ \phi_k^m (y):m \in [3^k ] \}$  
be a set of positive nondecreasing functions
on $\rr_+ = \{ y: y \ge 0 \}$.
Suppose that $\Phi_k$ satisfies a system $\sD$ of 
difference inequalities specified by a set of rooted labelled trees
$\{\sT_k^m: m \in [3^k ]\}$, 
such that all inequalities contain no advanced variables on their right side.
If the associated linear program $L_k^{\sD} ( \la )$ for  $\la > 1$
has a positive feasible solution with principal variables
$\{c_k^m\}$ then,
for all $m \in [3^k ]$, 
\beql{eq501}
\phi_k^m (y) \ge \Delta c_k^m \la^y ~,~~\mbox{for all}~~y \ge 0,
\eeq
with
\beql{eq502}
\Delta := \la^{-\nu} \frac{\min \{ \phi_k^m (0)\}}{\max\{ c_k^m\}}, 
\eeq
and $\nu$ is the largest backward time-shift of a variable in $\sD$.
\end{theorem}

\paragraph{Proof.}
Suppose that the set of functions $\Phi_k$ satisfies the system
$\sD := \{ D_k^m : m \in [3^k] \}$
of difference inequalities. Set
$$\mu := \min \{ \beta : \phi_k^{m'} (y- \beta )
~\mbox{appears on right side of some} ~ D_k^m \} 
$$
and 
$$
\nu := \max \{ \beta : \phi_k^{m'} (y- \beta ) 
~\mbox{appears on right side of some} ~ D_k^m \} ~.
$$
The hypothesis of no advanced variables in $\sD$ means that
$\mu > 0.$
Now the inequalities \eqn{eq501} hold for all $m \in [3^k]$, 
 on the initial interval $[0,\nu ]$, since
the definition of $\Delta$ gives
\beql{eq503}
\phi_k^m (y) \ge \phi_k^m (0) \ge \Delta \max \{c_k^m \} \lambda^\nu \ge 
\Delta c_k^m \la^y ~\mbox{for}~~ y \in [0, \bv],
\eeq
using the monotonicity and inequality properties of $\phi_k^m (y)$.

We now prove that \eqn{eq501} holds for all $m \in [3^k]$
on the interval 
$y \in [0, \nu + j \mu ]$ by induction on $j \ge 0$.
It holds for the base case $j=0$ by \eqn{eq503}.

For the induction step, suppose that \eqn{eq501} holds for $j$
and we are to prove it for $j+1$.
It suffices to  consider a given 
$y \in [\nu + j\mu, \nu + (j+1)\mu]$. The induction step
consists, schematically, of showing
\begin{eqnarray}
\phi_k^m(y) &\geq & \sum_{D_k^m(EL)} 
\mbox{nested-min}[ \phi_k^{m'}(y + \beta')]
~\label{eq504} \\
&\geq& \sum_{\sT_k^m(EL)} \mbox{nested-min}[c_k^{m'} \lambda^{y + \beta'}]
~\label{eq505} \\
&\geq& \Delta c_k^m \lambda^y  ~\label{eq506}.
\end{eqnarray}
Here \eqn{eq504} represents schematically the inequality $D_k^m(EL)$, 
with the right side actually being a nested series of minimizations.
Each function $\phi_k^{m'}(y + \beta')$ that appears on the
right side of \eqn{eq504} has $-\nu \le \beta' \le -\mu$, hence
$$ 0 \le j \mu \le y + \beta' \le \nu + j \mu,$$
so the induction hypothesis applies to each such  term.

The induction hypothesis gives
$$
\phi_k^{m'}(y + \beta') \ge \Delta c_k^{m'} \lambda^{y + \beta'}= 
\Delta \lambda^y(c_k^{m'}\lambda^{\beta'}).
$$
Substituting thes inequalities in \eqn{eq504} term by term yields the
right side of \eqn{eq505}, because
the nested minimization on the right side of \eqn{eq504} involves only the
operations of addition and minimization and these operations
are both monotone in each variable appearing in them; 
also the structure $\sT_k^{m}(EL)$ in \eqn{eq505} is the
tree structure of the inequality $D_k^m(EL)$. 
Now let $f(y)$ represent the value
of the right side of \eqn{eq505} as a function of $y$. Each minimization
on the right side of \eqn{eq505} corresponds to a m-vertex $\bv$ of 
$\sT_k^{m}(EL)$; we let $f_{\bv}(y)$ equal the value of this
minimization expression as a function of $y$. Next we
can apply the inequalities in $L_k^{EL} (\lambda)$ in a suitable
order to prove that
$$
f_{\bv}(y) \ge \Delta \lambda^y (a_{\bv}\lambda^{\beta(\bv)}).
$$
for all m-vertices; the order starts with the innermost 
minimization and works outward. At the last step we reach the
root vertex and obtain
$$
f(y) \ge \Delta \lambda^y c_k^m \lambda^{\beta(\bw_0)}= \Delta c_k^m \lambda^y,
$$
since $\beta(\bw_0)=0.$
This gives the right side of \eqn{eq506}. Since this
holds for all $k \in [3^m]$, this  completes the
induction step.
~~~$\bsq$

We now prove the main Theorem \ref{th21} by combining the results of \S3--\S5.

\paragraph{Proof of Theorem \ref{th21}.}
Theorem \ref{th32} shows that any set of positive 
nondecreasing functions $\Phi_k$ that satisfies
 the inequality system $\sI_k$ also satisfies the 
derived inequality system $\sI_k (EL)$ which has 
inequalities with no advanced variables on their right sides.
The family of linear programs associated to this 
inequality system in \S4 is denoted $L_k^{EL} (\lambda )$.

Suppose now that for a given $\la > 1$ the 
inequality system $L_k^{NT} (\la )$ has a 
feasible solution with principal variables
$\{c_k^m: m \in [3^k ] \}$.
Theorem \ref{th41} established that the linear program 
$L_k^{EL} (\la )$ has a positive feasible solution with the same 
principal variables and the same value of $\la$.

Theorem \ref{th51} then applies to the system $\sI_k (EL)$ 
to show that any positive feasible solution of 
$L_k^{EL} (\la )$ yields the bounds, for all $m \in [3^k]$,
\beql{eq507}
\phi_k^m (y) \ge \Delta c_m^k \la^y  ~\mbox{for all}~~y \ge 0,
\eeq
with
\beql{eq508}
\Delta := \la^{-\nu} \frac{ \min \{\phi_k^m (0) \}     }
{\max \{c_k^m\}} 
\eeq
where $\nu$ is the largest backwards timeshift.

For the system $\Phi_k = \{\phi_k^m (y): m \in [3^k] \}$ coming from the 
$3x+1$ problem, we have by (P1) that $\phi_k^m (0) \ge 1$.
We also have $\la \le 2$ and the maximum retarded term $\nu \le 2$.
Thus we have
$$\Delta \ge \Delta_1 = \frac{1}{4 \max{ \{ c_k^m \} } }
$$
which, with \eqn{eq507}, implies the desired bound \eqn{eq217}.~~~$\bsq$

\noindent\paragraph{Remark.} 
Theorem~\ref{th21} has the counterintuitive feature that
iterating the inequalities seems potentially to
strengthen, rather than weaken, the exponential lower
bound obtained. It allows the possibility
that the linear program $L_k^{EL}(\lambda)$ has a  positive
feasible solutions for a larger value of 
$\lambda$ than is obtainable using the original
linear program family $L_k^{NT}(\lambda)$.
However we believe this cannot occur, and that
the exponent obtained from $L_k^{NT}(\lambda)$ is
the largest possible for positive monotone solutions
to the original difference inequalities 
$\sI_k$. This is discussed at the end of \S6.

%
%
%
%
%
%
%
%

\section{$3X+1$ Lower Bounds}
\setcounter{equation}{0}

 We obtain lower bounds for the number $\pi_1(x)$ of integers below $x$
that eventually iterate to $1$ under the $3x+1$ function.

\begin{theorem}
For each positive $a \not\equiv 0 (\bmod~3)$ the function 
$$\pi_a(x) :=| \{ 1 \le n \le x: ~\mbox{Some}~~ T^{(j}(n) = a.\}|$$ 
satisfies, for all sufficiently large $x \ge x_0(a)$, 
$$ \pi_a(x) \ge x^{0.84}. $$   
\end{theorem}

\noindent\paragraph{Proof.} This follows from Theorem~\ref{th21}, by
finding a positive feasible solution by computer to the 
linear program family $L_{k}^{NT}(\lambda)$
for  $k=11$,
for $\lambda = 1.7922310,$ 
see Table 2 below. This  yields the exponent
$\gamma = \log_2 \lambda \approx 0.84175.$ ~~~$\bsq$.

Table 2
gives data on the 
 bounds for the optimal $\lambda$ for
$L_{k}^{NT}$ for $2 \le k \le 11$.
For $1\le k \le 9$ these are taken from \cite{AL95b};
the new values for $k=10, 11$ were computed by D. Applegate.

%
%
\begin{table}[hp]\centering
\begin{tabular}{|r|r|r|r|r|r|r|}
\hline
\multicolumn{1}{|c|}{$k$} & \multicolumn{1}{c|}{$\gamma_k$} &
\multicolumn{1}{c|}{$\lambda_k$} &
\multicolumn{1}{c|}{$C_k^{max}$} &
\multicolumn{1}{c|}{$\bar{c}_{k,k}$} &
\multicolumn{1}{c|}{$\bar{c}_{k-1,k}$} &
\multicolumn{1}{c|}{$\bar{c}_{k,k} - \bar{c}_{k-1,k}$}
\\ \hline
 2 & 0.4365880 & 1.3534010 & 1.8316920 & 1.5237640 & 1.0000000 & 0.5237640 \\
\hline
 3 & 0.6112620 & 1.5275960 & 3.4881908 & 2.1014900 & 1.6994294 & 0.4020606 \\
\hline
 4 & 0.6891080 & 1.6122870 & 5.4951954 & 2.7869040 & 2.4010985 & 0.3858055 \\
\hline
 5 & 0.7335790 & 1.6627590 & 9.0756176 & 3.4648343 & 3.0771822 & 0.3876521 \\
\hline
 6 & 0.7608180 & 1.6944520 & 12.8769418 & 3.9667005 & 3.5825321 & 0.3841684 \\
\hline
 7 & 0.7825670 & 1.7201900 & 20.1963763 & 4.8122983 & 4.4061650 & 0.4061333 \\
\hline
 8 & 0.8031960 & 1.7449630 & 29.1315157 & 5.2028179 & 4.8181536 & 0.3846643 \\
\hline
 9 & 0.8168300 & 1.7615320 & 43.3394210 & 5.8102043 & 5.4164870 & 0.3937173 \\
\hline
10 & 0.8295450 & 1.7771270 &  64.9801068 & 6.4567870 & 6.0648572 & 0.3919298 \\
\hline
11 & 0.8417560 & 1.7922310 & 98.4009647 & 7.1552344 & 6.7695583 & 0.3856761 \\
\hline
\end{tabular}
\caption{NLP  Lower Bounds: No truncation of advanced terms}~\label{table2}
\end{table}

The last three columns in Table~\ref{table2} give some average quantities
formulated in \cite{AL95b}. Define
$$
\bar{c}_{k, k} := \frac{1}{3^{k-1}} \sum_{m \in [3^k]} c_k^m
$$
and
$$
\bar{c}_{k-1,k} := \frac{1}{3^{k-2}}\sum_{m \in [3^{k-1}]} \bar{c}_{k-1}^m.
$$
Adding up all the inequalities in $L_k^{NT}(\lambda)$ leads to
$$
\bar{c}_{k, k} \le \lambda^{-2}\bar{c}_{k, k} + 
\frac{1}{3}(\lambda^{\alpha-1} + \lambda^{\alpha-2})\bar{c}_{k-1,k}.
$$
In \cite{AL95b} it was noted that
a  necessary and sufficient condition for a bound like
$\pi_1(x) > x^{1-\epsilon}$ to hold for 
each $\epsilon > 0$ and all sufficiently  large $x$
is that 
$\lambda_k \to 2$
as $k \to \infty$, and this in turn would follow from the
existence of feasible solutions with
$$
\frac{c_{k-1,k}}{c_{k,k}} \to 1~~\mbox{as}~~
k \to \infty.
$$
Table~\ref{table2} gives more empirical data on these quantities.

The supremum of the  exponential lower bounds that can be extracted from the
linear program family $L_k^{NT}(\lambda)$ is given by  
$\lambda_k$, the supremum of values of $\lambda$ for
which $L_k^{NT}(\lambda)$ has a feasible solution. These values
satisfy $\lambda_{k} \le \lambda_{k+1},$
because given a feasible solution to $L_k(\lambda)$  with
principal variables $c_k^m$ one can define  
$$c_{k+1}^{m+ j\cdot 3^k} := c_k^m ~~\mbox{ for}~~ 0 \le j \le 2,$$ 
and obtain a feasible solution to 
$L_{k+1}^{NT}(\lambda)$. It remains an  open problem
to show that the values $\lambda_k$ are strictly increasing
in $k$. As already noted in \cite{AL95b}, showing that $\lambda_k \to 2$
as $k \to \infty$ would imply a lower bound
$\pi_a (x) \ge x^{1 - \epsilon}$ 
holds for each positive $\epsilon$, for each $a \not\equiv 0~(\bmod~3)$
and all  sufficiently
large $x \ge x_0(a)$.

We  now relate the
linear program system $L_k^{NT}(\lambda)$ used here 
to the linear program system  denoted $L_{\lambda}^{NT}$
 in \cite{AL95b}. These 
 two linear program systems are equivalent
in the sense of Theorem~\ref{th21}; namely,
the set of $\lambda$
for which they have a strictly positive feasible solution coincide.
To see this, observe
first that if $L_k^{NT}(\lambda)$ has a feasible solution,
then it has a strictly positive feasible solution. One may have to
modify 
the auxiliary variables, which might be negative,
while holding the principal variables fixed.
However the auxiliary variables  can be forced to their maximal
values in terms of the principal variables without affecting feasiblity.
Such a feasible solution has  all  values at least $1$, so strict positivity
is attained, and this solution also 
satisfies $L_{\lambda}^{NT}$.
Conversely, given a positive feasible solution to $ L_{\lambda}^{NT}  $,
it can be multiplicatively rescaled to have objective
function value $c_1^2= 1$, and this gives  a feasible solution
to  $   L_k^{NT}(\lambda)$, on taking $C_k^{max} := \max{ \{c_k^m\}}.$

We conclude the paper by discussing the possibility that
the lower bound obtained in Theorem~\ref{th21} give the largest that
is implied by the difference inequalities $\sI_k$.
This would follow if one could exhibit a positive monotone solution
to $\sI_k$ that has a growth rate matching the lower bound.
Such a pure exponential lower bound could potentially be
constructed from a solution to $L_k^{NT}(\lambda_k)$.
Two conditions must hold: 

(1) The supremum $\lambda_k$ is attained. That is, $L_k^{NT}(\lambda_k)$
has a feasible solution.

(2) At the supremum value $\lambda_k$ , there exists a feasible solution
in which all of the principal inequalities (L1)-(L3) hold with
equality. 

If conditions (1), (2) hold, then the functions 
$\phi_k^m(y) = c_k^m \lambda_k^y$
would satisfy $\sI_k$ with equality for all times $y \ge 2$,
and would constitute a positive monotone solution to $\sI_k$
attaining the best lower bound given by Theorem~\ref{th21}.
Experimentally this is the case for $k \le 11$.

Regarding condition (1), $L_k^{NT}(\lambda_k)$ could fail to have a feasible
solution at the supremum value $\lambda_k$ only if  the objective function
value as $\lambda \to \lambda_k$ from below diverges to $\infty$, so
some variables $c_k^m$ become unbounded. 
Regarding condition (2), the complementary slackness conditions for
an optimal solution would for a  generic linear program
of this type force all the principal inequalities (L1)-(L3) to hold
with equality. This would happen, for example, if 
there were an optimal solution at which all variables
$c_k^m$ took distinct values. We think it likely that properties
(1), (2) hold for all $k \ge 2$, but this may be difficult to prove.

The supremum linear
program  $L_k^{NT}(\lambda_k)$ has a finite
optimal objective function value $C_k^{max}$ 
provided that condition (1) holds, as we now assume:
denote this value  by $\tilde{C}_k^{max}$.
The value $\tilde{C}_k^{max}$ 
has an interesting meaning: it measures the minimal spread 
attainable in the values
of $c_k^m$, while normalizing these variables by $\min{ \{c_k^m\}}=1$.
This quantity shows up in the  constant $\Delta_1$ in Theorem~\ref{th21}.
One may view the  value $\tilde{C}_k^{max}$ as a quantitative
measure of a rate of  ``mixing'' 
between congruence classes $(\bmod~3^k)$ that the $3x+1$ function
produces.
The fourth column of Table 2
indicates that the quantity  $\tilde{C}_k^{max}$ 
exists for $k \le 11$, and it appears to grow
exponentially with $k$.

{\bf Acknowledgments.}  The authors thank 
David Applegate for computations given in Tables 1 and 2.

%
%
%
%
%
%
%
%

\section*{Appendix:  Inequalities for $k=2$}

The  case $k=2$ is the only case where the derived inequalities
$\sI_k(EM)$ and the linear program family $L_k^{EM}(\lambda)$
can be easily written down. There are three functions
$\Phi_2 := \{ \phi_2^2 (y), \phi_2^5 (y), \phi_2^8 (y): y \ge 0 \} ~.$
Recall that $\alpha = \log_2 3 \approx 1.585$.
The inequalities $\sI_2$ are
\begin{eqnarray*}
\phi_2^2 (y) & \ge &
\phi_2^8 (y-2) + \min [\phi_2^2 (y+ \alpha -2), \phi_2^5 (y+ \alpha-2 ),
\phi_2^8 (y + \alpha -2 )] \,, \\
\phi_2^5 (y) & \ge & \phi_2^2 (y-2) ~, \\
\phi_2^8 (y) & \ge & \phi_2^5 (y-2) + \min [\phi_2^2 (y+\alpha -1), 
\phi_2^5 (y+\alpha - 1),\phi_2^8 (y+\alpha -1) ] \,.
\end{eqnarray*}
Of these, only the  inequality for $\phi_2^8 (y)$
contains advanced terms on its right side.

%
%

\begin{figure}[htb]
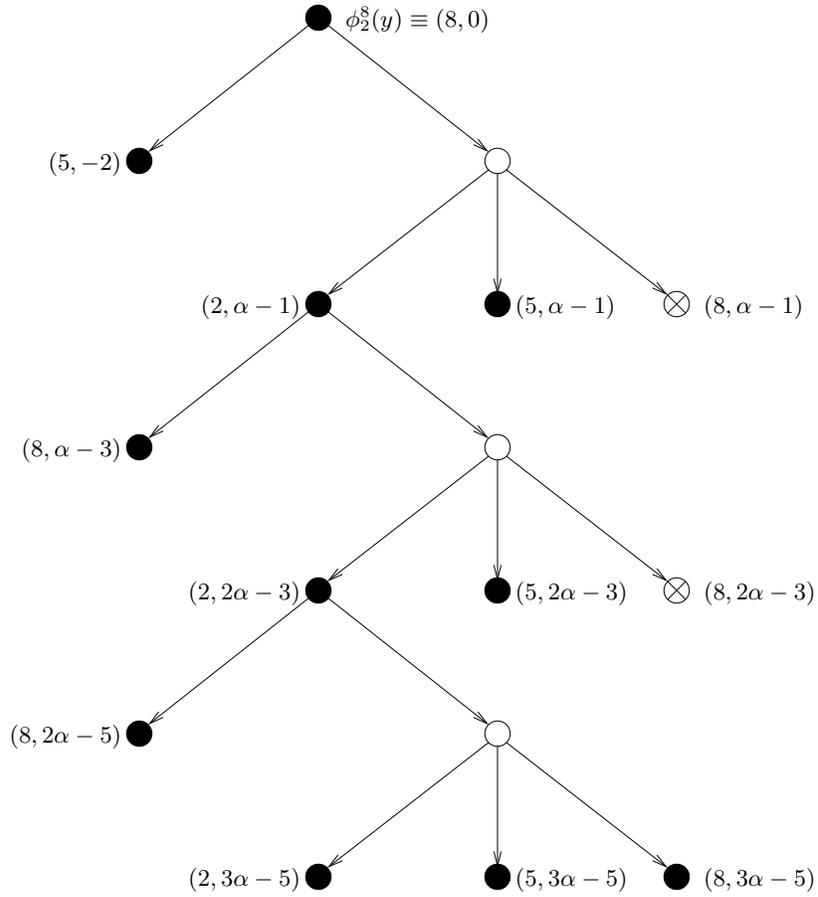

\begin{center}
\input kafgA1.pstex_t
\end{center}
\caption{Tree $\sT_2^8 (EL)$ (Nodes  marked $\otimes$ are deleted nodes)}
\end{figure}~\label{figure3}

The corresponding inequality $\sI_2^8(EL)$ 
has three leaves of nested minimization.
The corresponding tree $\sT_2^8 (EL)$ is pictured in
Figure 3, 
with the deleted nodes indicated.
The tree $\sT_2^8 (EL )$ has three  $m$-nodes and eight leaf nodes.
We let  $a_1$, $a_2$, $a_3$ be the auxiliary variables for the leaf nodes,
numbered as in Figure 3, and its associated inequalities are:
$$\phi_2^8 (y) \ge \phi_2^5 (y-2) + 
\min [\phi_2^8 (y+\alpha -3) +M_1 (y), \phi_2^2 (y+\alpha -3) ] \,,
$$
in which
$$M_1 (y): = \min [\phi_2^8 (y+2 \alpha -5 ) + M_2 (y) , 
\phi_2^5 (y+2\alpha -5 )] \,,
$$
and
$$M_2 (y) := \min [\phi_2^2 (y+3 \alpha -5 ), 
\phi_2^5 (y+ 3 \alpha -5 ) , \phi_2^8 (y+3 \alpha-5 )] \,.
$$

The inequalities in the 
linear program $L_2^{EL} (\lambda )$ for
the three trees $\sT_2^m(EL)$ with $m = 2, 5$ and $8$ 
are given in Table~\ref{table4}; they are 
associated to the leaves of these trees,
identified by their labels in Table~\ref{table4}.

\begin{table}[htb]
$$
\begin{array}{|l|l|l|} \hline
\mbox{Tree} & \mbox{Leaf node label} & 
\multicolumn{1}{c}{\mbox{Inequality}} \\ \hline
~ & ~ \\
~& (8, - 2) & c_2^2 \le c_2^8 \la^{-2} + 
a_1^{'}\la^{\alpha - 2} \\
~ & ~ & ~ \\
\sT_2^2(EL) & (2,~ \alpha - 2) & 
a_1^{'}\la^{\alpha - 2}  \le  c_2^2\la^{\alpha - 2} \\
~ & ~ & ~ \\
~ & (5,~\alpha - 2) & a_1^{'} \la^{\alpha - 2} \le  c_2^5\la^{\alpha - 2} \\
~ &~ & ~ \\
~ & (8,~\alpha - 2) &  a_1^{'} \la^{\alpha - 2} \le c_2^8\la^{\alpha - 2} \\
~ &~ & ~ \\
\hline

~ &~ & ~ \\
\sT_2^5(EL)    & (2,~-2) & c_2^5 \le c_2^2 \la^{-2} \\
~ &~ & ~ \\
\hline

~ &~ & ~ \\
~ &  (5,~ -2) & c_2^8 \le c_2^5 \la^{-2} + a_1 \la^{\alpha -1} \\
~ &~ & ~ \\
~ & (8,~\alpha -3 ) & a_1 \la^{\alpha -1} \le 
c_2^8 \la^{\alpha -3} + a_2 \la^{2 \alpha -3} \\
~ &~ & ~ \\
~ & (2,~ \alpha -3) & a_1 \la^{\alpha -1} \le 
c_2^2 \la^{\alpha -3} \\
~ &~ & ~ \\
\sT_2^8(EL)   & (8, ~2\alpha -5) & a_2 \la^{2 \la -3} \le 
c_2^8 \la^{2 \alpha -5} + a_3 \la^{3 \alpha -5} \\
~ &~ & ~ \\
~ & (2, ~2 \alpha -5) & a_2 \la^{2 \alpha -3} \le 
c_2^2 \la^{2 \alpha -5} \\ 
~ & ~ \\
~ & (2, ~3 \alpha -5) & a_3 \la^{3 \alpha -5} \le 
c_2^2 \la^{3 \alpha -5} \\ 
~ & ~ \\
~ & (5, ~3 \alpha -5) & a_3 \la^{3 \alpha -5} \le 
c_2^5 \la^{3 \alpha -5} \\ 
~ & ~ \\
~ & (8, ~3 \alpha -5) & a_3 \la^{3 \alpha -5} \le 
c_2^8 \la^{3 \alpha -5} \\ 
\hline
\end{array}
$$
\caption{$L_2^{EL}(\la )$ inequalities for trees
$\sT_2^m (EL)$.}~\label{table4}
\end{table}

{\tt
\begin{tabular}{ll}
email: &Ilia.Krasikov@brunel.ac.uk \\
&jcl@research.att.com \\
\end{tabular}
}

\end{document}